\documentclass[12pt]{article}
\usepackage{amssymb,amsthm,amsmath,latexsym}
\usepackage{url}
\newtheorem{thm}{Theorem}
\newtheorem{prop}[thm]{Proposition}

\theoremstyle{remark}
\newtheorem{rem}[thm]{Remark}
\newcommand{\FF}{\mathbb{F}}

\newcommand{\allone}{\mathbf{1}}
\newcommand{\cC}{\mathcal{C}}
\newcommand{\cD}{\mathcal{D}}

\DeclareMathOperator{\wt}{wt}

\DeclareMathOperator{\supp}{supp}

\begin{document}

\title{New extremal singly even self-dual codes of
lengths $64$ and $66$}

\author{
Damyan Anev\thanks{
Faculty of Mathematics and Informatics,
Konstantin Preslavski University of Shumen,
Shumen, 9712, Bulgaria.},
Masaaki Harada\thanks{
Research Center for Pure and Applied Mathematics,
Graduate School of Information Sciences,
Tohoku University, Sendai 980--8579, Japan.}
and
Nikolay Yankov\thanks{
Faculty of Mathematics and Informatics,
Konstantin Preslavski University of Shumen,
Shumen, 9712, Bulgaria.}
}

\maketitle

\begin{abstract}
For lengths $64$ and $66$, 
we construct extremal singly even self-dual 
codes with weight enumerators for which no extremal
singly even self-dual codes were previously known to exist.
We also construct new $40$ inequivalent
extremal doubly even self-dual $[64,32,12]$ codes
with covering radius $12$ meeting the Delsarte bound.
\end{abstract}

\section{Introduction}

A (binary) $[n,k]$ {\em code} $C$ is a $k$-dimensional vector subspace
of $\FF_2^n$,
where $\FF_2$ denotes the finite field of order $2$.
All codes in this note are binary.
The parameter $n$ is called the {\em length} of $C$.
The {\em weight} $\wt(x)$ of a vector $x$ is
the number of non-zero components of $x$.
A vector of $C$ is a {\em codeword} of $C$.
The minimum non-zero weight of all codewords in $C$ is called
the {\em minimum weight} of $C$. An $[n,k]$ code with minimum
weight $d$ is called an $[n,k,d]$ code.
The \textit{dual code} $C^{\perp}$ of a code
$C$ of length $n$ is defined as
$
C^{\perp}=
\{x \in \FF_2^n \mid x \cdot y = 0 \text{ for all } y \in C\},
$
where $x \cdot y$ is the standard inner product.
A code $C$ is called \textit{self-dual} if $C = C^{\perp}$.
A self-dual code $C$ is {\em doubly even} if all
codewords of $C$ have weight divisible by four, and 
{\em singly even} if there is at least one codeword $x$
with $\wt(x) \equiv 2 \pmod 4$.
It is  known that a self-dual code of length $n$ exists
if and only if  $n$ is even, and
a doubly even self-dual code of length $n$
exists if and only if $n$ is divisible by $8$.

Let $C$ be a singly even self-dual code.
Let $C_0$ denote the subcode of $C$ consisting of 
codewords $x$ with $\wt(x) \equiv 0 \pmod 4$.
The {\em shadow} $S$ of $C$ is defined to be $C_0^\perp \setminus C$.
Shadows for self-dual codes were introduced by Conway and
Sloane~\cite{C-S}
in order to give the largest possible minimum weight 
among singly even self-dual codes, and to provide
restrictions on the weight enumerators of 
singly even self-dual codes.
The largest possible minimum weights
among singly even self-dual codes of length $n$ were given
for $n \le 72$ in~\cite{C-S}.
The possible weight enumerators of singly 
even self-dual codes with
the largest possible minimum weights 
were given 
in~\cite{C-S} and~\cite{Dougherty1997} for $n \le 72$.
It is a fundamental problem to find which weight enumerators 
actually occur
for the possible weight enumerators (see~\cite{C-S}).
By considering the shadows, 
Rains~\cite{Rains} showed that
the minimum weight $d$ of a self-dual code of length $n$
is bounded by
$d  \le 4 \lfloor{\frac {n}{24}} \rfloor + 6$
if $n \equiv 22 \pmod {24}$,
$d  \le 4  \lfloor{\frac {n}{24}} \rfloor + 4$
otherwise.
A self-dual code meeting the bound is called  {\em extremal}.

The aim of this note is to
construct extremal singly even self-dual 
codes with weight enumerators for which no extremal
singly even self-dual codes were previously known to exist.
More precisely, we construct
extremal singly even self-dual $[64, 32, 12]$ codes with
weight enumerators 
$W_{64,1}$ for $\beta=35$, and 
$W_{64,2}$ for $\beta \in \{19, 34, 42, 45,50\}$
(see Section~\ref{sec:WE} for $W_{64,1}$ and $W_{64,2}$).
These codes are constructed as self-dual neighbors of
extremal four-circulant singly even self-dual codes.
We construct
extremal singly even self-dual $[66, 33, 12]$ codes with
weight enumerators
$W_{66,1}$ for $\beta \in \{7,58,70,91,93\}$, and
$W_{66,3}$ for $\beta \in \{22,23\}$
(see Section~\ref{sec:WE} for $W_{66,1}$ and $W_{66,3}$).
These codes are constructed from extremal singly even
self-dual $[64,32,12]$ codes by the method given 
in~\cite{Tsai92}.
We also demonstrate that
there are at least $44$ inequivalent
extremal doubly even self-dual $[64,32,12]$ codes
with covering radius $12$ meeting the Delsarte bound.

All computer calculations in this note
were done with the help of
the algebra software {\sc Magma}~\cite{Magma}
and 
the computer system Q-extensions~\cite{Q}.

\section{Weight enumerators of extremal singly even self-dual
codes of lengths 64 and 66}
\label{sec:WE} 

The possible weight enumerators $W_{64,i}$ and $S_{64,i}$
of extremal singly even self-dual $[64,32,12]$ codes
and their shadows are given in~\cite{C-S}:
\begin{align*}
&
\begin{cases}
W_{64,1} =& 1+(1312 + 16 \beta)y^{12} +(22016 -64 \beta) y^{14}
+ \cdots,\\
S_{64,1} =& y^4+(\beta-14)y^8 + (3419-12\beta)y^{12}+ \cdots,
\end{cases}\\
&
\begin{cases}
W_{64,2} =& 1+(1312 + 16 \beta) y^{12}+(23040 - 64 \beta) y^{14}
+ \cdots,\\
S_{64,2} =& \beta y^8 + (3328-12\beta)y^{12} + \cdots,
\end{cases}
\end{align*}
where $\beta$ are integers with
$14 \le \beta \le 104$ for $W_{64,1}$ and 
$0 \le \beta \le 277$ for $W_{64,2}$.
Extremal singly even self-dual codes with weight enumerator $W_{64,1}$ 
are known for
\[
\beta \in 
\left\{\begin{array}{l}
14, 16, 18, 20, 22, 24, 25, 26, 28, 29, 30, 32, \\
34, 36, 38, 39, 44, 46, 53, 59, 60, 64, 74 
\end{array}\right\}
\]
(see~\cite{CHK64}, \cite{Kay17b}, \cite{KYP16} and \cite{Yan14}). 
Extremal singly even self-dual codes with weight enumerator $W_{64,2}$ 
are known for
\[
\beta \in 
\left\{\begin{array}{l}
0,1,\ldots,41,44, 48, 51, 52, 56, 58, 64, 65, 72, \\
80, 88, 96, 104, 108, 112, 114, 118, 120, 184
\end{array}\right\}
\setminus\{19, 31, 34, 39\}
\]
(see~\cite{CHK64}, \cite{Kay17b}, \cite{Yan14} and \cite{YIL17}).

The possible weight enumerators $W_{66,i}$ and $S_{66,i}$ of
extremal singly even self-dual $[66,33,12]$ codes and their shadows are
given in~\cite{Dougherty1997}:
\begin{align*}
&\left\{\begin{array}{l}
W_{66,1} = 1 + (858 + 8\beta)y^{12} + (18678 - 24\beta)y^{14} 
+ \cdots,\\
S_{66,1} = \beta
 y^{9}+(10032-12\beta)y^{13} 
+\cdots,
\end{array}
\right.\\
&\left\{\begin{array}{l}
W_{66,2} = 1 + 1690y^{12} + 7990y^{14} 
+\cdots,\\
S_{66,2} = y + 9680y^{13} 
+\cdots,\\
\end{array}
\right.\\
&\left\{\begin{array}{l}
W_{66,3}= 1 + (858 + 8\beta)y^{12} + (18166 - 24\beta)y^{14} 
+ \cdots,\\
S_{66,3}= y^{5}+(\beta-14)y^{9}+(10123-12\beta)y^{13}+\cdots,
\end{array}\right.
\end{align*}
where $\beta$ are integers with 
$0\leq\beta\leq 778$ for $W_{66,1}$ and
$14\leq\beta\leq 756$ for $W_{66,3}$.
Extremal singly even self-dual codes with weight enumerator $W_{66,1}$
are known for
\[
\beta \in \{0,1,\ldots,92,94, 100, 101, 115\} \setminus \{4, 7, 
58, 70, 91\}
\]
(see~\cite{CKO}, \cite{HNY07}, \cite{Kay17b}, \cite{YLGI15} 
and \cite{YIL17}).
Extremal singly even self-dual codes with weight enumerator $W_{66,2}$
are known (see~\cite{HNY07} and \cite{Tsa99}).
Extremal singly even self-dual codes with weight enumerator $W_{66,3}$
are known for
\[
\beta\in \{24,25,\ldots,92\} \setminus \{65,
68, 69, 72, 89, 91\}
\]
(see~\cite{KY13}, \cite{Kay17b}, \cite{KYP16} and \cite{KYS15}).

\section{Extremal four-circulant
singly even self-dual $[64,32,12]$ codes}\label{sec:64S}

An $n \times n$ circulant matrix has the following form:
\[
\left(
\begin{array}{ccccc}
r_0&r_1&r_2& \cdots &r_{n-1} \\
r_{n-1}&r_0&r_1& \cdots &r_{n-2} \\
\vdots &\vdots & \vdots && \vdots\\
r_1&r_2&r_3& \cdots&r_0
\end{array}
\right),
\]
so that each successive row is a cyclic shift of the previous one.
Let $A$ and $B$ be $n \times n$ circulant matrices.
Let $C$ be a $[4n,2n]$ code with generator matrix of the following form:
\begin{equation} \label{eq:GM}
\left(
\begin{array}{ccc@{}c}
\quad & {\Large I_{2n}} & \quad &
\begin{array}{cc}
A & B \\
B^T & A^T
\end{array}
\end{array}
\right),
\end{equation}
where $I_n$ denotes the identity matrix of order $n$
and $A^T$ denotes the transpose of $A$.
It is easy to see that $C$ is self-dual if
$AA^T+BB^T=I_n$.
The codes with generator matrices of the form~\eqref{eq:GM}
are called {\em four-circulant}.

Two codes are {\em equivalent} if one can be obtained from the other by a
permutation of coordinates.
In this section, we give a classification of extremal
four-circulant singly even
self-dual $[64,32,12]$ codes.
Our exhaustive search found all distinct extremal four-circulant
singly even self-dual $[64,32,12]$ codes,
which must be checked further for equivalence to
complete the classification.
This was done by considering all pairs of
$16 \times 16$ circulant matrices $A$ and $B$ satisfying
the condition that
$AA^T+BB^T=I_{16}$, the sum of
the weights of the first rows of $A$ and $B$
is congruent to $1 \pmod 4$
and the sum of the weights is greater than or equal to $13$.
Since a cyclic shift of the first rows gives an equivalent code,
we may assume without loss of generality that
the last entry of the first row of $B$ is $1$.
Then our computer search shows that
the above distinct extremal four-circulant
singly even self-dual $[64,32,12]$ codes are divided into
$67$ inequivalent codes.

\begin{prop}
Up to equivalence, there are $67$ extremal four-circulant
singly even self-dual $[64,32,12]$ codes.
\end{prop}

We denote the $67$ codes by $C_{64,i}$ $(i=1,2,\ldots,67)$.
For the $67$ codes $C_{64,i}$,
the first rows $r_A$ (resp.\ $r_B$) of the circulant matrices $A$
(resp.\ $B$)
in generator matrices~\eqref{eq:GM} are listed in Table~\ref{Tab:64S}.
We verified that the codes $C_{64,i}$ have
weight enumerator $W_{64,2}$, where $\beta$ are also listed
in Table~\ref{Tab:64S}.

\begin{table}[thbp]
\caption{Extremal four-circulant singly even self-dual $[64,32,12]$ codes }
\label{Tab:64S}
\begin{center}
{\scriptsize
\begin{tabular}{c|c|c|c}
\noalign{\hrule height0.8pt}
Codes & $r_A$&$r_B$ & $\beta$\\
\hline
$C_{64,1}$&(0000001100111111)&(0001011010101111)&0\\
$C_{64,2}$&(0000010101111101)&(0010011010111011)&0\\
$C_{64,3}$&(0000011001101111)&(0010110101011011)&0\\
$C_{64,4}$&(0000000001011111)&(0001001100101011)&8\\
$C_{64,5}$&(0000000010101111)&(0011011011110111)&8\\
$C_{64,6}$&(0000000011010111)&(0000100110011011)&8\\
$C_{64,7}$&(0000000011010111)&(0000101100010111)&8\\
$C_{64,8}$&(0000000011010111)&(0011101110101111)&8\\
$C_{64,9}$&(0000000110111111)&(0101101111111111)&8\\
$C_{64,10}$&(0000001001011101)&(0001000101011011)&8\\
$C_{64,11}$&(0000001100011111)&(0010101011011111)&8\\
$C_{64,12}$&(0000001100011111)&(0010111011011011)&8\\
$C_{64,13}$&(0000001100111011)&(0001101011101111)&8\\
$C_{64,14}$&(0000001101111111)&(0011101111011111)&8\\
$C_{64,15}$&(0000010000111101)&(0010111011011111)&8\\
$C_{64,16}$&(0000010001011111)&(0001110101101111)&8\\
$C_{64,17}$&(0000010110111011)&(0001101110001111)&8\\
$C_{64,18}$&(0000000100011111)&(0010111111110011)&16\\
$C_{64,19}$&(0000000100111101)&(0000101011000111)&16\\
$C_{64,20}$&(0000000110010111)&(0001001111111111)&16\\
$C_{64,21}$&(0000000111001111)&(0010101110111101)&16\\
$C_{64,22}$&(0000000111001111)&(0010110110111011)&16\\
$C_{64,23}$&(0000001000101111)&(0011101011110111)&16\\
$C_{64,24}$&(0000001011100011)&(0010101111110111)&16\\
$C_{64,25}$&(0000001011100011)&(0011011011111011)&16\\
$C_{64,26}$&(0000010010011111)&(0010110011101111)&16\\
$C_{64,27}$&(0000011001101111)&(0001001011011111)&16\\
$C_{64,28}$&(0000011011011111)&(0010010101011101)&16\\
$C_{64,29}$&(0000011011100111)&(0001011111001011)&16\\
$C_{64,30}$&(0000011101111111)&(0101101110110111)&16\\
$C_{64,31}$&(0000101110111111)&(0011101011110111)&16\\
$C_{64,32}$&(0000000000100111)&(0001011101101011)&24\\
$C_{64,33}$&(0000000001011011)&(0010010101101011)&24\\
$C_{64,34}$&(0000000100111111)&(0001001000101011)&24\\
$C_{64,35}$&(0000000101001011)&(0010010110011011)&24\\
$C_{64,36}$&(0000000101001011)&(0010011001011011)&24\\
$C_{64,37}$&(0000000110111111)&(0000001000100111)&24\\
$C_{64,38}$&(0000001001111111)&(0010101111001011)&24\\
$C_{64,39}$&(0000001100011111)&(0001010011111111)&24\\
$C_{64,40}$&(0000001100011111)&(0001110011110111)&24\\
$C_{64,41}$&(0000010001011111)&(0010101111001111)&24\\
$C_{64,42}$&(0000010001101111)&(0011001110101111)&24\\
$C_{64,43}$&(0000010011101111)&(0001011101100111)&24\\
$C_{64,44}$&(0000010101010111)&(0001010111101111)&24\\
$C_{64,45}$&(0000010101010111)&(0010110011111011)&24\\
$C_{64,46}$&(0000010101110111)&(0000101111110011)&24\\
$C_{64,47}$&(0000010101110111)&(0001011101101011)&24\\
$C_{64,48}$&(0000011011110111)&(0101101110111111)&24\\
$C_{64,49}$&(0000000001001011)&(0000111010110111)&32\\
$C_{64,50}$&(0000000001100111)&(0001001111100011)&32\\
\noalign{\hrule height0.8pt}
\end{tabular}
}
\end{center}
\end{table}

\setcounter{table}{0}
\begin{table}[thb]
\caption{Extremal four-circulant singly even self-dual $[64,32,12]$
 codes (continued)}
\begin{center}
{\scriptsize
\begin{tabular}{c|c|c|c}
\noalign{\hrule height0.8pt}
Codes & $r_A$&$r_B$ & $\beta$\\
\hline
$C_{64,51}$&(0000001010111011)&(0001011111100111)&32\\
$C_{64,52}$&(0000010101011111)&(0001101111000111)&32\\
$C_{64,53}$&(0000010101111101)&(0010110010110111)&32\\
$C_{64,54}$&(0000011010111111)&(0000101110011101)&32\\
$C_{64,55}$&(0000101011101011)&(0001011111001011)&32\\
$C_{64,56}$&(0000000000100111)&(0001011010111011)&40\\
$C_{64,57}$&(0000000010101101)&(0001001011011011)&40\\
$C_{64,58}$&(0000001000011101)&(0000100101111011)&40\\
$C_{64,59}$&(0000001110011111)&(0001010111101101)&40\\
$C_{64,60}$&(0000011000111111)&(0001010111101101)&40\\
$C_{64,61}$&(0000011011001111)&(0000101010111111)&40\\
$C_{64,62}$&(0000100111011111)&(0001010101011011)&40\\
$C_{64,63}$&(0000001001101011)&(0001010011001101)&48\\
$C_{64,64}$&(0000000001011011)&(0001011000101111)&56\\
$C_{64,65}$&(0000010111011111)&(0010100101011011)&56\\
$C_{64,66}$&(0000101110011101)&(0001000101111111)&64\\
$C_{64,67}$&(0000000001011111)&(0001011111110111)&72\\
\noalign{\hrule height0.8pt}
\end{tabular}
}
\end{center}
\end{table}

\section{Extremal self-dual $[64,32,12]$ neighbors of $C_{64,i}$}
\label{sec:SEd12}

Two self-dual codes $C$ and $C'$ of length $n$
are said to be {\em neighbors} if $\dim(C \cap C')=n/2-1$.
Any self-dual code of length $n$ can be reached
from any other by taking successive neighbors (see~\cite{C-S}).
Since every self-dual code $C$ of length $n$ contains the all-one vector
$\allone$, $C$ has $2^{n/2-1}-1$ subcodes $D$ of codimension $1$
containing $\allone$. Since $\dim(D^\perp/D)=2$, there are two
self-dual codes rather than $C$ lying between $D^\perp$ and $D$.
If $C$ is a singly even
self-dual code of length divisible by $8$,
then $C$ has two doubly even self-dual neighbors
(see~\cite{BP}).
In this section, we construct extremal self-dual
$[64,32,12]$ codes by considering self-dual neighbors.

For $i=1,2,\ldots,67$,
we found all distinct extremal singly even self-dual neighbors of 
$C_{64,i}$, which are equivalent to none of the $67$ codes.
Then we verified that these codes are divided into
$385$ inequivalent codes $D_{64,i}$ $(i=1,2,\ldots,385)$.
These codes $D_{64,i}$ are constructed as
\[
\langle (C_{64,j} \cap \langle x \rangle^\perp), x \rangle.
\]
To save space, 
the values $j$, the supports $\supp(x)$ of $x$,
the values $(k,\beta)$ in the weight enumerators $W_{64,k}$
are listed in \\
``\url{http://www.math.is.tohoku.ac.jp/~mharada/Paper/64-SE-d12.txt}''
\\
for the $385$ codes.
For extremal singly even self-dual $[64,32,12]$
codes with weight enumerators for which no extremal
singly even self-dual codes were previously known to exist,
$j$, $\supp(x)$ and $(k,\beta)$ are list in Table~\ref{Tab:64nei}.
Hence, we have the following:

\begin{prop}
There is an extremal singly even self-dual $[64, 32, 12]$ code with
weight enumerator 
$W_{64,1}$ for $\beta=35$, and 
$W_{64,2}$ for $\beta \in \{19, 34, 42, 45,50\}$.
\end{prop}

\begin{table}[thbp]
\caption{Extremal singly even self-dual $[64,32,12]$ neighbors}
\label{Tab:64nei}
\begin{center}
{\footnotesize
\begin{tabular}{c|c|l|c}
\noalign{\hrule height0.8pt}
Codes & $j$ & \multicolumn{1}{c|}{$\supp(x)$} & $(k,\beta)$\\
\hline
$D_{64,138}$&24&$\{1,2,3,38,42,43,45,46,48,54,56,57\}$ &$(2,19)$\\
$D_{64,270}$&49&$\{1,2,8,32,38,41,48,49,50,53,55,61\}$ &$(1,35)$\\
$D_{64,283}$&52&$\{1,2,4,33,36,37,41,43,46,51,61,64\}$ &$(2,42)$\\
$D_{64,293}$&56&$\{3,7,9,10,11,37,43,53,57,58,62,64\}$ &$(2,34)$\\
$D_{64,314}$&64&$\{6,8,26,37,38,40,43,46,48,59,61,63\}$&$(2,50)$\\
$D_{64,329}$&65&$\{1,6,8,9,37,47,50,52,57,60,63,64\}$  &$(2,45)$\\
\hline
$D_{64,  1}$& 1&$\{4,7,9,34,38,40,45,46,47,50,51,53\}$&$  (2,0)$\\
$D_{64,  2}$& 1&$\{3,37,38,47,48,50,52,53,54,59,60,63\}$&$(2,0)$\\
$D_{64, 12}$& 4&$\{2,4,5,16,17,38,40,46,56,57,60,62\}$&$  (2,0)$\\
$D_{64, 19}$& 4&$\{2,3,6,7,9,35,41,49,55,56,57,61\}$&$    (2,0)$\\
$D_{64, 22}$& 4&$\{2,33,34,35,38,39,42,45,48,52,61,62\}$&$(2,0)$\\
$D_{64, 33}$& 6&$\{8,9,10,16,17,33,44,45,54,55,59,61\}$&$ (2,0)$\\
$D_{64, 44}$& 6&$\{1,3,6,33,36,38,39,45,47,55,57,59\}$&$  (2,0)$\\
$D_{64, 58}$& 8&$\{1,3,5,16,17,35,36,38,42,44,54,59\}$&$  (2,0)$\\
$D_{64, 66}$& 8&$\{4,6,9,34,36,39,41,42,48,51,57,63\}$&$  (2,0)$\\
$D_{64, 68}$& 8&$\{3,6,9,33,36,37,38,49,56,57,60,62\}$&$  (2,0)$\\
$D_{64, 84}$&13&$\{1,4,5,35,37,38,41,44,53,60,61,62\}$&$  (2,0)$\\
$D_{64, 95}$&13&$\{2,4,9,34,35,40,42,47,49,52,59,64\}$&$  (2,0)$\\
$D_{64,108}$&15&$\{2,16,17,37,43,48,49,52,54,57,58,64\}$&$(2,0)$\\
$D_{64,115}$&16&$\{1,3,6,7,8,41,45,46,49,50,57,60\}$&$    (2,0)$\\
$D_{64,136}$&21&$\{3,16,17,33,34,37,42,44,47,51,52,56\}$&$(2,0)$\\
$D_{64,143}$&26&$\{1,2,9,34,37,38,41,48,57,58,59,64\}$&$  (2,0)$\\
$D_{64,191}$&35&$\{1,2,6,8,10,33,37,46,54,59,60,63\}$&$   (2,0)$\\
$D_{64,240}$&47&$\{2,4,7,9,13,16,17,44,56,59,62,64\}$&$   (2,0)$\\
$D_{64,254}$&48&$\{1,2,5,7,8,35,36,37,45,47,49,63\}$&$    (2,0)$\\
\hline
$D_{64,14}$& 4&$\{1,7,8,35,36,37,41,43,46,49,51,53\}$&$(1,14)$\\
$D_{64,383}$ &67&$\{1,33,34,36,37,38,40,41,47,49,50,53,55,59,61,63\}$&$(2,40)$\\
\noalign{\hrule height0.8pt}
\end{tabular}
}
\end{center}
\end{table}

Now we consider the extremal doubly even self-dual neighbors
of $C_{64,i}$ $(i=1,2,3)$.
Since the shadow has minimum weight $12$,
the two doubly even self-dual neighbors 
$\cC^1_{64,i}$ and  $\cC^2_{64,i}$
are extremal doubly even self-dual $[64,32,12]$ codes
with covering radius $12$  (see~\cite{CHK64}).
Thus, six extremal doubly even self-dual $[64,32,12]$ codes with
covering radius $12$ are constructed.
In addition, 
among the $385$ codes $D_{64,i}$ $(i=1,2,\ldots,385)$, 
the $19$ extremal singly even self-dual codes $D_{64,j}$ 
have shadow of minimum weight $12$, where
\[
j \in \{
  1,  2, 12, 19, 22, 33, 44, 58, 66, 68, 84, 95,
108,115,136,143,191,240,254\}.
\]
The constructions of the $19$ codes $D_{64,j}$ are listed 
in Table~\ref{Tab:64nei}.
Their two doubly even self-dual neighbors 
$\cD_{64,j}^1$ and $\cD_{64,j}^2$ 
are extremal doubly even self-dual $[64,32,12]$ codes
with covering radius $12$.
We verified that there are
the following equivalent codes among 
the four codes in~\cite{CHK64},
the six codes $\cC_{64,i}^1$, $\cC_{64,i}^2$
and the $38$ codes $\cD_{64,j}^1$, $\cD_{64,j}^2$, where
\[
\cD_{64, 22}^2 \cong \cD_{64,68}^2,
\cD_{64, 33}^2 \cong \cD_{64,84}^2,
\cD_{64, 44}^2 \cong \cD_{64,95}^2,
\cD_{64,136}^2 \cong \cD_{64,143}^2,
\]
where $C \cong D$ means that $C$ and $D$ are equivalent,
and there is no other pair of equivalent codes.
Therefore, we have the following proposition.

\begin{prop}
There are at least $44$ inequivalent
extremal doubly even self-dual $[64,32,12]$ codes
with covering radius $12$ meeting the Delsarte bound.
\end{prop}

In order to distinguish two doubly even neighbors
$\cD_{64,i}^1$ and $\cD_{64,i}^2$ $(i=68,84,95,143)$,
we list in Table~\ref{Tab:64neiDE} the supports $\supp(x)$ for 
the $8$ codes, where
$\cD_{64,i}^1$ and $\cD_{64,i}^2$ 
are constructed as
$\langle (D_{64,i} \cap \langle x \rangle^\perp), x \rangle$.

\begin{table}[thbp]
\caption{Extremal doubly even self-dual $[64,32,12]$ neighbors}
\label{Tab:64neiDE}
\begin{center}
{\footnotesize
\begin{tabular}{c|l}
\noalign{\hrule height0.8pt}
Codes & \multicolumn{1}{c}{$\supp(x)$} \\
\hline
$\cD_{64, 68}^1$&$\{1,4,7,34,35,36,47,54,55,58,60,63\}$\\
$\cD_{64, 68}^2$&$\{1,4,5,6,30,42,45,47,54,56,58,64\}$\\
$\cD_{64, 84}^1$&$\{16,17,33,39,43,46,48,49,51,54,58,64\}$\\
$\cD_{64, 84}^2$&$\{1,2,6,33,35,38,40,42,52,57,59,60\}$\\
$\cD_{64, 95}^1$&$\{1,2,6,33,35,38,40,42,52,57,59,60\}$\\
$\cD_{64, 95}^2$&$\{3,33,38,41,45,47,51,53,58,60,62,64\}$\\
$\cD_{64,143}^1$&$\{1,4,10,40,43,46,52,54,58,61,62,63\}$\\
$\cD_{64,143}^2$&$\{1,31,34,42,44,45,46,50,51,52,54,62\}$\\
\noalign{\hrule height0.8pt}
\end{tabular}
}
\end{center}
\end{table}

\section{Four-circulant singly even
self-dual $[64,32,10]$ codes and self-dual neighbors}
\label{sec:SEd10}

Using an approach similar to that given in Section~\ref{sec:64S},
our exhaustive search found all distinct
four-circulant singly even self-dual $[64,32,10]$ codes.
Then our computer search shows that
the distinct four-circulant singly even
self-dual $[64,32,10]$ codes are divided into
$224$ inequivalent codes.

\begin{prop}
Up to equivalence, there are $224$
four-circulant singly even self-dual $[64,32,10]$ codes.
\end{prop}

We denote the $224$ codes by $E_{64,i}$ $(i=1,2,\ldots,224)$.
For the codes,
the first rows $r_A$ (resp.\ $r_B$) of the circulant matrices
$A$ (resp.\ $B$)
in generator matrices~\eqref{eq:GM} can be obtained from\\
``\url{http://www.math.is.tohoku.ac.jp/~mharada/Paper/64-4cir-d10.txt}''.

The following method for constructing self-dual neighbors
was given in~\cite{CHK64}.
For $C=E_{64,i}$ $(i=1,2,\ldots,224)$,
let $M$ be a matrix  whose rows are the codewords of weight $10$ in $C$.
Suppose that there is a vector $x$ of even weight such that
\begin{equation}\label{Eq:N}
M x^T = \text{\boldmath$1$}^T.
\end{equation}
Then $C^0=\langle x\rangle^\perp\cap C$
is a subcode of index $2$ in $C$.
We have self-dual neighbors $\langle C^0,x\rangle$ and $\langle
C^0,x+y\rangle$ of $C$
for some vector $y\in C\setminus C^0$,
which have no codeword of weight $10$ in $C$.
When $C$ has a self-dual neighbor $C'$ with minimum weight $12$,
there is a vector $x$ satisfying~\eqref{Eq:N}
and we can obtain $C'$ in this way.
For $i=1,2,\ldots,224$,
we verified that there is a unique vector satisfying~\eqref{Eq:N}
and $C$ has two self-dual neighbors, where
$C^0$ is a doubly even $[64,31,12]$ code.
In this case, the two neighbors are automatically doubly even.
Hence, we have the following:

\begin{prop}
There is no extremal singly even self-dual $[64,32,12]$ neighbor 
of $E_{64,i}$ for $i=1,2,\ldots,224$.
\end{prop}

\section{Extremal singly even self-dual $[66, 33, 12]$ codes}

The following method for constructing singly even self-dual
codes was given in~\cite{Tsai92}.
Let $C$ be a self-dual code of length $n$.
Let $x$ be a vector of odd weight.
Let $C^0$ denote the subcode of $C$
consisting of all codewords which are orthogonal to $x$.
Then there are cosets
$C^1,C^2,C^3$ of $C^0$ such that ${C^0}^\perp = C^0 \cup C^1 \cup
C^2 \cup C^3$, where $C = C^0  \cup C^2$ and $x+C = C^1 \cup C^3$.
It was shown in~\cite{Tsai92} that
\begin{equation}\label{eq:code}
C(x)= (0,0,C^0) \cup  (1,1,C^2) \cup  (1,0,C^1) \cup  (0,1,C^3)
\end{equation}
is a self-dual code of length $n+2$.
In this section, we construct new extremal singly even self-dual codes
of length $66$
using this construction from the extremal singly even self-dual $[64,
32, 12]$ codes
obtained in Sections~\ref{sec:64S} and \ref{sec:SEd12}.

Our exhaustive search shows that
there are $1166$ inequivalent extremal singly even
self-dual $[66, 33, 12]$ codes constructed as the codes $C(x)$ 
in~\eqref{eq:code}
from the codes $C_{64,i}$ $(i=1,2,\dots, 67)$.
$1157$ codes of the $1166$ codes have weight enumerator $W_{66,1}$
for
$\beta\in \{7,8,\ldots,92\} \setminus \{9, 11\}$,
$3$ of them have weight enumerator $W_{66,3}$ for $\beta \in \{30, 49, 54\}$,
and
$6$ of them have weight enumerator  $W_{66,2}$.
Extremal singly even self-dual $[66, 33, 12]$ codes with weight
enumerator
$W_{66,1}$ for $\beta \in \{7, 58, 70,91\}$ are constructed for the
first time.
For the four weight enumerators $W$, 
as an example, codes $C_{66,i}$ with weight enumerators $W$ are given
$(i=1,2,3,4)$.
We list in Table~\ref{Tab_66new}
the values $\beta$ in $W$,
the codes $C$ and the vectors $x=(x_1,x_2,\ldots,x_{32})$ 
of $C(x)$ in~\eqref{eq:code}, where $x_j=1$ $(j=33,\ldots,64)$.

\begin{table}[thb]
\caption{Extremal singly even self-dual $[66, 33, 12]$ codes}
\label{Tab_66new}
\begin{center}
{\footnotesize
\begin{tabular}{c|c|c|c|c}
\noalign{\hrule height0.8pt}
Codes & $\beta$&$W$&$C$&$(x_1,\dots, x_{32})$\\
\hline
$C_{66,1}$&  7&$W_{66,1}$&$C_{64,1}$&(01101101101010010111111010101100)\\
$C_{66,2}$& 58&$W_{66,1}$&$C_{64,56}$&(00001101100000011000110000011100)\\
$C_{66,3}$& 70&$W_{66,1}$&$C_{64,66}$&(00100110011011001001011100000010)\\
$C_{66,4}$& 91&$W_{66,1}$&$C_{64,67}$&(00001110110111110000011101000010)\\
\hline
$D_{66,1}$&22 &$W_{66,3}$ &$D_{64,14}$ &(10100011100100110111101010011111)\\
$D_{66,2}$&23 &$W_{66,3}$ &$D_{64,14}$ &(10111100111100000100101000100011)\\
$D_{66,3}$&93 &$W_{66,1}$
	 &$D_{64,383}$&(10100101011110010011001101001101)\\
\noalign{\hrule height0.8pt}
\end{tabular}
}
\end{center}
\end{table}

By applying the construction given in~\eqref{eq:code} to $D_{64,i}$,
we found more extremal singly even self-dual $[66,33,12]$ 
codes $D_{66,j}$ with weight enumerators for which no extremal
singly even self-dual codes were previously known to exist.
For the codes $D_{66,j}$, 
we list in Table~\ref{Tab_66new}
the values $\beta$ in the weight enumerators $W$,
the codes $C$ and the vectors $x=(x_1,x_2,\ldots,x_{32})$ 
of $C(x)$ in~\eqref{eq:code}, where $x_i=1$ $(i=33,\ldots,64)$.
Hence, we have the following:

\begin{prop}
There is an extremal singly even self-dual $[66, 33, 12]$ code with
weight enumerator
$W_{66,1}$ for $\beta \in \{7,58,70,91,93\}$,
and weight enumerator
$W_{66,3}$ for $\beta \in \{22,23\}$.
\end{prop}

\begin{rem}
The code $D_{66,1}$ has the smallest value $\beta$
among known extremal singly even self-dual $[66,33,12]$ codes
with weight enumerator $W_{66,3}$.
\end{rem}

\bigskip
\noindent
{\bf Acknowledgment.}
This work was supported by JSPS KAKENHI Grant Number 15H03633.



\end{document}